\documentclass[12pt]{article}

\usepackage{bussproofs}
\EnableBpAbbreviations

\usepackage{url}

\def\HDS{\vrule width0pt height2.3ex depth1.05ex\displaystyle}
\def\f#1#2{{{\HDS #1}\over{\HDS #2}}}

\begin{document}

\title{{\LARGE G\" odel's Notre Dame Course}}
\author{{\sc Milo\v s Ad\v zi\' c} and {\sc Kosta Do\v sen}
\\[1ex]
{\small Faculty of Philosophy, University of Belgrade}\\[-.5ex]
{\small \v Cika Ljubina 18-20, 11000 Belgrade, Serbia, and}\\
{\small Mathematical Institute, Serbian Academy of Sciences and Arts}\\[-.5ex]
{\small Knez Mihailova 36, p.f.\ 367, 11001 Belgrade, Serbia}\\[.5ex]
{\small email: milos.adzic@gmail.com, kosta@mi.sanu.ac.rs}}
\date{\small April 2016}
\maketitle

\begin{abstract}
\noindent This is a companion to a paper by the authors entitled ``G\" odel's natural deduction'', which presented and made comments about the natural deduction system in G\" odel's unpublished notes for the elementary logic course he gave at the University of Notre Dame in 1939. In that earlier paper, which was itself a companion to a paper that examined the links between some philosophical views ascribed to G\" odel and general proof theory, one can find a brief summary of G\" odel's notes for the Notre Dame course. In order to put the earlier paper in proper perspective, a more complete summary of these interesting notes, with comments concerning them, is given here.
\end{abstract}

\noindent {\small \emph{Keywords:} propositional logic, predicate logic}

\vspace{1ex}

\noindent {\small \emph{Mathematics Subject Classification
(2010):} 01A60 (History of mathematics and mathematicians, 20th century), 03-01 (Mathematical logic and foundations, instructional exposition), 03-03 (Mathematical logic and foundations, historical)}

\vspace{1ex}

\noindent {\small \emph{Acknowledgements.} Work on this paper was
supported by the Ministry of Education, Science and Technological Development of Serbia.
We are very grateful indeed to Gabriella Crocco for making some of G\" odel's
unpublished notes available to us. The decipherment and publication of some of these notes
is part of the project {\it Kurt G\" odel Philosopher: From
Logic to Cosmology}, which is directed by her and funded by the French
National Research Agency (project ANR-09-BLAN-0313). She and Antonio Piccolomini d'Aragona were also very kind to invite us to the workshop {\it Inferences and Proofs}, in Marseille, in May 2016, where the second-mentioned of us delivered a talk based partly on this paper, and enjoyed their especial hospitality. We would like to thank Professor Patricia Blanchette for her careful editorial work.
We are grateful to the Institute for Advanced Study in Princeton for granting us, through the office of its
librarian Mrs Marcia Tucker and its archivist Mr Casey Westerman, the permission to
use G\" odel's unpublished notes; we were asked to give credit for that with the following text:
``All works of Kurt G\" odel used with permission. Unpublished Copyright
(1934-1978) Institute for Advanced Study. All rights reserved by Institute
for Advanced Study.''}

\vspace{1.5ex}

\section*{Introduction}
This paper is a companion to the paper \cite{DA16a}, which is itself a companion to the paper \cite{DA16}, where we examined links between some philosophical views ascribed to G\" odel and general proof theory. In \cite{DA16a} we presented and made comments about the natural deduction system in G\" odel's unpublished notes for the main elementary logic course he gave in his career.

Among G\" odel's unpublished writings at the Princeton University Library, two sets of notes for elementary logic courses have been preserved (see Section 1.II of \cite{Daw05}, pp.\ 152-154). He gave the first of these courses at the University of Vienna in the summer of 1935 (see \cite{Daw97}, p.\ 108). According to \cite{Daw05} (p.\ 153), the notes for the Vienna course are about: ``\ldots truth tables, predicate logic, Skolem normal form, the Skolem-L\" owenheim theorem, the decision problem, and set theory''.

The second of these courses of G\" odel's was for graduate students at the University of Notre Dame in the spring semester of 1939 (see \cite{Daw97}, pp.\ 135-136, and \cite{Daw}). Our aim here is to summarize the notes for this second course, and make some comments concerning them. We gave a brief summary of these notes in \cite{DA16a}, and very brief summaries of them and a few fragments may be found in Dawson's writings mentioned above and in \cite{CN09}, where three bigger extracts are published. We believe that these notes deserve a more detailed treatment, which we want to supply here. We hope that thereby the earlier papers that this paper accompanies will be put in proper perspective. We will not however mention absolutely all that G\" odel's notes contain. To produce a readable and relatively short summary, we supply only samples of his material, and select matters we find important, rather than make a complete inventory.

\section*{Major problems and branches of logic}
G\" odel's logic course at the University of Notre Dame was elementary, and he was initially somewhat skeptical about his competence to give such a course. In a letter to Karl Menger (see \cite{Go03a}, pp.\ 117-119), he gave the following reasons for his doubts: ``As concerns the program of lectures, I think that I am presently not very well suited [for] giving an elementary course of lectures, on account of insufficient knowledge of English, insufficient experience at elementary lectures and insufficient time for preparation.''

In spite of these doubts, G\" odel prepared a short introduction to logic very well thought out, covering quite a lot of ground and friendly towards beginners. In just one semester, he introduced his students to propositional logic, predicate logic and also briefly at the end to the theory of types.

In many important respects, G\" odel's course, as much of his work in logic, was influenced by \cite{HA28}. This textbook, first published in 1928, was based on a series of lectures that Hilbert gave at the University of G\" ottingen between 1917 and 1922. The second edition of this book, published in 1938, contains, besides smaller additions, G\" odel's completeness proof for first-order predicate logic from his doctoral thesis of 1929.

In the notes for his course, G\" odel intended to introduce his students to three major problems that logicians were to deal with in the XXth century, and through these problems to nearly all the major branches of logic, which were all shaped through his own fundamental results. He dealt in the notes with the completeness problem through a completeness proof for propositional logic and mentioned the completeness of first-order predicate logic. Thereby he made a first, modest, step in the direction of model theory. G\" odel stressed the importance of completeness, and how the successful treatment of it in modern logic makes a tremendous advance in comparison with what we had in the old logical tradition. He dealt also with another kind of completeness---the functional completeness of sets of connectives.

G\" odel dealt in the notes with the decidability problem for tautologies through truth tables, and referred briefly and picturesquely to Turing machines (without mentioning them by name). Thereby he made a first step in the direction of recursion theory, which was founded in the 1930s with his decisive contribution.

He dealt in the notes with the independence problem, i.e.\ with the independence of axioms, for a system of the Hilbert type for propositional logic. Thereby he made in a certain sense a first step in the direction of set theory, where this problem was central. He made a step towards set theory also in the notes for the last part of the course, where he dealt with the paradoxes and the theory of types.

G\" odel worked on independence in set theory at the time he gave the Notre Dame course.
Together with this elementary course, he was giving at the same time at Notre Dame a more advanced course on his recent results in set theory about the consistency of the Axiom of Choice and the Continuum Hypothesis with Zermelo-Fraenkel set theory. Although consistency is mentioned here, these results of G\" odel may be phrased as showing that the negations of these two postulates are independent from Zermelo-Fraenkel set theory. The independence of the two postulates was shown by Paul Cohen in the 1960s, while G\" odel claimed (in a letter to Wolfgang Rautenberg; see \cite{Go03a}, pp.\ 179-183) that he had an uncompleted proof for the independence of the Axiom of Choice, on which he might have worked during his time in Notre Dame.

Perhaps for lack of time, G\" odel in the notes for his course did not envisage dealing with the consistency problem, which is central for Hilbertian proof theory, though G\" odel's most famous incompleteness results are so important for that problem and that theory. However, with his natural deduction system based on sequents, G\" odel made in the notes for the course a first small step in the direction of proof theory as well. This direction may be understood as being rather that of general proof theory than that of Hilbertian proof theory.

The four problems we mentioned---completeness, decidability, independence and consistency---together with the division of logic into the branches where the study of each of them, respectively, dominates or characterizes the whole branch---model theory, recursion theory, set theory and proof theory---have become standard for education in logic in the second half of the XXth century. G\" odel's course may today look pretty standard (in accordance with \cite{TASL}), but, at the time when he gave it, we believe that it was still a novelty to present logic in this manner. (And it was a feat to do it in a single semester.) These problems may be found in a nutshell already in \cite{HA28}, where consistency, independence and completeness are singled out in Section I.12, while Section III.12 is about decidability. G\" odel was however, as we indicated, one of the main figures, if not the main figure, to give modern logic its profile by dealing so successfully with these problems.

\section*{The contents of the course}
We refer the reader to the papers \cite{Daw05} (Section 1.II, pp.\ 153-154) and \cite{CN09} (beginning of Section~6, p.\ 77) for a description of the state of G\" odel's unpublished notes in English for the Notre Dame course, which we are now going to summarize, and comment upon a little bit. We give below a rough table of contents for these notes, which we made without taking account of everything in them. We did not attempt to classify in our table, nor will we summarize below, lists of exercises, examples, repetitions, sketchy notes, the text involving the Gabelsberger shorthand, which, according to \cite{Daw97} (p.\ 136), contains examination questions, and notes apparently foreign to the course (dealing mainly with religious matters). The division of the course into sections, which we numerate with the section sign \S $\,$, is ours, as well as their titles, with terms that are not necessarily G\" odel's:
\begin{tabbing}
\hspace{2em}\=\textbf{Propositional logic}\\*[.5ex]
\>\S 1.\hspace{.8em}\=Failure of traditional logic\\
\>\S 2.\>Connectives\\
\>\S 3.\>Tautologies and decidability\\
\>\S 4.\>Axiom system for propositional logic\\
\>\S 5.\>Functional completeness\\
\>\S 6.\>Completeness for propositional logic\\
\>\S 7.\>Independence of the axioms\\
\>\S 8.\>Comments on completeness\\
\>\S 9.\>Sequents and natural deduction system\\[1ex]

\>\textbf{Predicate logic}\\*[.5ex]
\>\S 10.\>First-order languages and valid formulae\\
\>\S 11.\>Decidability and completeness in predicate logic\\
\>\S 12.\>Axiom system for predicate logic\\
\>\S 13.\>Remarks on the term ``tautology''\\
\>\S 14.\>``Thinking machines''\\
\>\S 15.\>Existential presuppositions\\
\>\S 16.\>Classes and relations\\
\>\S 17.\>Paradoxes and type theory
\end{tabbing}

Our summary with comments in the remainder of this paper is divided into subsections numerated by reference to this table. In this summary, as we announced at the end of the Introduction, we don't cover everything G\" odel has in the course, but only what appears to us to be prominent, original or interesting.

\section*{Propositional logic}
\textbf{\S 1. Failure of traditional logic.} G\" odel's notes for his lectures on logic begin by considering how inadequate it was before it appeared as a mathematical subject. It failed concerning completeness, and decidability too. It gave instead ``a more or less arbitrary selection from the infinity of the laws of logic''.  Occasional claims that everything can be deduced from this or that law have never been proved or even clearly formulated. It is a great achievement of modern logic, according to G\" odel, that the completeness problem was precisely formulated, after which it became possible to prove systems of propositional and predicate logic complete.

G\" odel finds that in traditional logic, with all its failures, the contribution of the Stoics to propositional logic is more fundamental than the Aristotelian contribution with syllogistic figures and moods. (The same opinion may be found in \S 35 of the famous scholarly work on this matter \cite{L50}.) After stressing the fundamental importance of understanding logical form, contrary to Aristotle, as it is understood in propositional logic, G\" odel introduces the symbols, taken over from Russell (and as in \cite{WR10}), for the usual logical connectives: negation $\sim$, conjunction $.$ (as in $p\, .\: q$), disjunction $\vee$, implication $\supset$ and equivalence $\equiv$. He remarks that in traditional logic disjunction has been understood in the exclusive manner. After explaining how complex propositions are built with connectives, G\" odel considers briefly the Polish prefix notation for propositional formulae. He remarks that the connective \emph{if} (without the accompanying \emph{then}) is used in ordinary language in the Polish way.

\vspace{2ex}

\noindent \textbf{\S 2. Connectives.} G\" odel then proceeds to discuss the meaning of the logical connectives, introducing truth tables and the principle of truth-functionality. After explaining material implication with care, he envisages an intensional implication tied to deduction, which he calls ``strict implication''. G\" odel discusses the difference between intensional and extensional, i.e.\ truth-functional or material, implication, and finds that the former would be worth studying later in the course if time permits.

He does mention at one place ``the logic of modalities'', and also uses the fish-hook symbol for strict implication $\prec$ (which in his handwriting does not differ much from $<$), but it is not excluded that he also envisaged something else. He could have had in mind intuitionistic implication, which by a translation equivalent with his translation of \cite{Go33a} corresponds to the S4 strict implication. Moreover, in \cite{G35}, G\" odel could have found a clear connection between deduction and intuitionistic implication.

If G\" odel did not envisage intuitionistic implication as his intensional implication tied to deduction, he should have done so, because, as later developments in general proof theory were to show, this implication is tied to deduction in a deeper way than strict implication. We have in mind here the characterization of intuitionistic implication through adjunction in categorial proof theory (which originates in \cite{Law69}; some references for this and related matters may be found in Section~5 of \cite{DA16}).

A very interesting side remark that G\" odel  makes about material implication in this context is that it corresponds as closely to \emph{if then} as a precise, mathematical, notion can correspond to an imprecise notion of ordinary language (see \cite{CN09}, p.\ 83, fn 14).

Then G\" odel discusses the tautologies $q\supset(p\supset q)$ and $\sim p\supset(p\supset q)$, and says that their apparent paradoxicality arises only from the intensional point of view. He discusses also the tautology $(p\supset q)\vee(q\supset p)$, which seems as paradoxical. The first two tautologies are kept as theorems in intuitionistic logic, in which implication is not truth-functional (as G\" odel found in \cite{Go32}), and is more intensional. Intuitionistic propositional logic may be understood as arising by introducing this more intensional implication, tied to deduction, and this is one of the main features of intuitionism, if not the main one. The third tautology, $(p\supset q)\vee(q\supset p)$, is, however, rejected in intuitionistic logic. This formula, which serves to axiomatize Dummett's intermediate logic of \cite{D59}, holds in G\" odel's intermediate logics of \cite{Go32}.
(Extract 1 in \cite{CN09}, which is transcribed from Notebook~I of G\" odel's notes, and represents the second version of something covered too in Notebook~0, covers \S\S 1-2.)

\vspace{2ex}

\noindent \textbf{\S 3. Tautologies and decidability.} Before speaking of \emph{tautologies}, G\" odel called them universally true or logically true formulae of propositional logic. He defines them as formulae that are true whatever propositions we put for their propositional letters. Then G\" odel considers a number of important tautologies and comments upon them.

G\" odel envisages verifying whether a propositional formula is a tautology in the standard manner, by going through all the possibilities in a truth table, and he explains why such a verification must have a truth table of $2^n$ lines for a formula with $n$ different variables. He notes however that in practice the number of cases having to be separately considered is smaller, because several cases may be dealt with in the same manner, and combined into one. He concludes that this simple solution of decidability for propositional logic is chiefly due to having only extensional connectives.

G\" odel notes that with strict implication decidability would have been much more complicated, and that it has been solved only recently ``under certain assumptions about strict \dots'' (here the text breaks, and a dozen pages are missing). It is not clear what G\" odel had in mind at this place. We suggested above that by ``strict implication'' G\" odel should have referred to intuitionistic implication. A decision procedure for the propositional modal system S5 could be found at that time in \cite{W33}, but if our suggestion is accepted, G\"odel could have had in mind here Gentzen's solution of the decidability problem for intuitionistic propositional logic from \cite{G35} (Section IV.1.2).

\vspace{2ex}

\noindent \textbf{\S 4. Axiom system for propositional logic.} G\" odel then proceeds to set up a formal system of the Hilbert type for propositional logic. As primitive connectives, he has disjunction and negation. As rules of inference, he has modus ponens, the rule of substitution, and the rules of defined symbols, which permit the replacement of a defined symbol by its definiens and vice versa. As axioms he has the following four axioms, which are the propositional axioms of \textit{Principia Mathematica} \cite{WR10} (Section I.A.$\ast$1) with one axiom omitted:
\begin{tabbing}
\hspace{11em}\= $p \supset (p \vee q)$,\\*[.5ex]
\>$(p \vee p) \supset p$,\\[.5ex]
\>$(p \vee q) \supset (q \vee p)$,\\[.5ex]
\>$(p \supset q) \supset ((r \vee p) \supset (r \vee q))$.
\end{tabbing}
(We don't follow G\" odel in omitting parentheses.) The axiom $(p\vee (q \vee r))\supset (q \vee (p \vee r))$ of \textit{Principia} is omitted because, as G\" odel remarks, Bernays showed that it is not independent from the other axioms (see \cite{Be26}; Bernays' proof is from his second habilitation thesis of 1918, about which see \cite{Z99}). The same axiom system for propositional logic as G\" odel's is in \cite{HA28} (Section I.10).

G\" odel verifies by discussing possible cases, and not by writing entire truth tables, that the four axioms are tautologies, and while introducing the rules of inference he established that by passing from premises to conclusions they preserve the property of being a tautology. Thus he has established the soundness of this system of propositional logic, i.e.\ that every theorem of the system is a tautology, while introducing the system. He proceeds to prove a number of basic theorems and to obtain many derived rules of inference.

As an interesting aside concerning axiom systems (which in the notes is displaced to occur later, within \S 5), G\" odel remarks that one might ask whether we really need both axioms and rules of inference, since implicational axioms themselves can already be taken to suggest rules. Rejecting this idea as entirely wrong, he points out the difference between material implication and deduction, echoing the moral of Lewis Carroll's Achilles and the Tortoise of \cite{Ca95}, which should be that there can be no deduction without rules of deduction. He also notes that even if it could be said that axioms suggest rules of inference, they do not uniquely determine these rules, since the same axiom could suggest different rules. For example, the axiom ${p \supset (p \vee q)}$ could suggest both of the following rules:
\[
\f{p}{p \vee q}\quad\quad\quad\f{\sim(p\vee q)}{\sim p}
\]

\vspace{2ex}

\noindent \textbf{\S 5. Functional completeness.} After that, G\" odel turns to the question of functional completeness. He remarks that instead of the set ${\{\vee, \sim\}}$ of primitive logical connectives, various other sets would do as well. He briefly mentions the Sheffer stroke, after which he proceeds to prove that the set ${\{\equiv, \sim\}}$ is not functionally complete. G\" odel gives a rather detailed proof of this result, carefully illustrating the important points by examples. He concludes, as a corollary, that exclusive disjunction, which is negated equivalence, and negation do not make a functionally complete set of connectives either.

Turning now to the proof that the set ${\{\vee, \sim\}}$ is functionally complete, instead of giving the proof in full generality for $n$-ary connectives, G\" odel chose to deal only with ternary truth-functional connectives. This is enough to get the idea of the proof, which is based on reading a disjunctive normal form of the formula $f(p,q,r)$, for $f$ an arbitrary ternary truth-functional connective, out of its truth table. Since the formula in disjunctive normal form contains only $\vee$, $.$ and $\sim$, and $.$ is definable in terms of $\vee$ and $\sim$, one can infer definability in terms of $\vee$ and $\sim$ alone.

\vspace{2ex}

\noindent \textbf{\S 6. Completeness for propositional logic.} When he introduced his system for propositional logic of the Hilbert type, G\" odel established at the same time, as we noted, the soundness of this system. He turns now towards the proof of the converse implication, i.e.\ completeness for this system. He presents this completeness proof with great detail, breaking it up into a series of lemmata, and providing also occasionally examples, to help the students understand the key points. The proof is based on what is sometimes called Kalmar's Lemma (see \emph{Hilfssatz}~3 in \cite{Ka35}, Section~8).

\vspace{2ex}

\noindent \textbf{\S 7. Independence of the axioms.} After the proof of the completeness theorem, G\" odel takes up the question of the independence of logical axioms. Using Bernays' method with matrices, he proves the second of the four axioms independent from the rest, noting the same could be done for the remaining three axioms. There is also a brief mention of the procedure of reducing a propositional formula to disjunctive normal form and conjunctive normal form.

\vspace{2ex}

\noindent \textbf{\S 8. Comments on completeness.} G\" odel had a favourable opinion about Gentzen's presentation of logic with sequents introduced in \cite{G35}. After proving completeness for the axiomatization of propositional logic of the Hilbert type, which we considered in the preceding section, G\" odel made some comments on the importance of this proof (considered in Section~3 of \cite{DA16a}, with a longer quotation from G\" odel's notes). It achieves something that traditional logic did not even conceive properly.

G\" odel says that it is not very important that this particular system has been proved complete. The idea is important. He says that he chose the system above because it has become standard, but other formal systems may also be proved complete, and among them some are \emph{aesthetically} more satisfactory than this system of the Hilbert type.

We have commented on the importance of this aesthetic point of view for G\" odel and mathematics in Section~2 of \cite{DA16} and Section~3 of \cite{DA16a}. We suggested that this is not a slight matter, and accords well with G\" odel's platonism.

In an aesthetically more satisfactory system, it will not happen that the very simple $p\supset p$ is proved from axioms more complicated than it. Such a system is a natural deduction system, which we are now going to present.

\vspace{2ex}

\noindent \textbf{\S 9. Sequents and natural deduction system.} G\" odel's natural deduction system in the notes for the Notre Dame course can be briefly described as being Ja\' skowski's system of \cite{J34} presented with Gentzen's sequents of \cite{G35}. G\" odel does not use the terms \emph{natural deduction} and \emph{sequent}; he mentions Gentzen, but not Ja\' skowski. This natural deduction system is presented in more detail, with comments, in Sections 4-6 of \cite{DA16a}, which are here summarized.

G\" odel's propositional language has the propositional letters $p,q,r,\ldots$, and the capital letters $P,Q,R,\ldots$ are schematic letters for the propositional formulae of the language with the primitive connectives of implication $\supset$ and negation $\sim$. The capital Greek letters $\Delta,\Gamma,\ldots$ are schematic letters for sequences of these propositional formulae. In a sequent $\Delta\rightarrow P$, which G\"odel conceives as a word in a language, the sequence of propositional formulae $\Delta$ should be finite.

Then G\" odel presents the axioms and rules of inference for his natural deduction system. As axioms he has only the identity law, i.e.\ all sequents of the form $P\rightarrow P$, and as rules of inference he has thinning:
\[
\f{\Delta\rightarrow Q}{P,\Delta\rightarrow Q}\quad\quad\quad\f{\Delta\rightarrow Q}{\Delta, P\rightarrow Q}
\]
the rules for introducing and eliminating implication:
\[
\f{\Delta, P\rightarrow Q}{\Delta\rightarrow P\supset Q}\quad\quad\quad\f{\Delta\rightarrow P\quad\quad\quad\Delta\rightarrow P\supset Q}{\Delta\rightarrow Q}
\]
and, finally, the rule of \emph{reductio ad absurdum}, in its strong, nonconstructive, version:
\[
\f{\Delta,\sim P\rightarrow Q\quad\quad\quad\Delta,\sim P\rightarrow \;\sim Q}{\Delta\rightarrow P}
\]

Then G\" odel claims that this system is complete, by which he means presumably that if $P$ is a tautology, then the sequent $\rightarrow P$, with the empty left-hand side, is provable. If we want to introduce other connectives, G\" odel's way of dealing with that would be to have rules corresponding to the definitions of these connectives in terms of $\supset$ and $\sim$. (He had an analogous approach with his systems of the Hilbert type; see \S 4 above.) At the end of the first part of the course, devoted to propositional logic, G\" odel says that he is sorry that he has ``no time left to go into more details about this Gentzen system''.

\section*{Predicate logic}
\noindent \textbf{\S 10. First-order languages and valid formulae.} G\" odel's treatment of first-order predicate logic is not short, but covers less major results than his treatment of propositional logic. In a motivating discussion of predicates and their role in first-order languages, he mentions the importance of predicates of arity greater than one, which were neglected in the old Aristotelian logic, and says that they are more important for the applications of logic in mathematics and elsewhere. After that he introduces quantifiers and offers some simple examples of formalization in a language of predicate logic, which may contain equality. For his axiomatic system, he takes the universal quantifier, written $(x)$, as primitive, defining the existential quantifier $(\exists x)$ in terms of it and negation in the usual way. He considers the equivalences obtained by permuting quantifiers of the same kind and the failure of equivalence when the existential quantifier is permuted with a universal one. The notions of free and bound variables, together with the notion of the scope of a quantifier, are introduced, and some important equivalences involving quantifiers, like renaming of bound variables and, later, Herbrand's laws of passage, are given.

\vspace{2ex}

\noindent \textbf{\S 11. Decidability and completeness in predicate logic.} The completeness of first-order predicate logic, from G\" odel's doctoral thesis, is only mentioned, but not proved. The same holds for the undecidability of predicate logic, and the decidability of monadic predicate logic (i.e.\ first-order predicate logic with only unary predicates).

\vspace{2ex}

\noindent \textbf{\S 12. Axiom system for predicate logic.} In setting up the formal system for predicate logic, G\" odel adds the axiom $(x) \varphi(x) \supset \varphi(y)$ to the propositional axioms, he extends the propositional rule of substitution to cover also substitution in predicate logic (which is pretty involved, because it involves predicate and propositional variables besides individual variables), and he adds the following rule for the universal quantifier:
\[
\f{\psi\supset\varphi(x)}{\psi\supset (x)\varphi(x)}
\]
with the proviso that the variable $x$ does not occur free in the formula $\psi$. On all these points, G\" odel is closely following the exposition in the book \cite{HA28} (Section III.5). A minor difference is that Hilbert and Ackermann take both quantifiers as primitive, and have also the dual axiom and rule of inference for the existential quantifier. (According to Hilbert and Ackermann, these axioms and rules are due to Bernays.) After that, G\" odel also considers a number of derived rules of inference.

\vspace{2ex}

\noindent \textbf{\S 13. Remarks on the term ``tautology''.} G\" odel uses the term ``tautology'' not only for the universally true formulae, i.e.\ logical truths, of propositional logic, but also for such formulae of predicate logic, which nowadays are called rather \emph{valid formulae}. (As for tautologies of propositional logic, he defines valid formulae apparently in a rather syntactical way, by appealing to substitution, and not by explicitly mentioning model theoretical interpretations.) This usage might perhaps be considered better than what has prevailed, because of its uniformity, but G\" odel says that using ``tautology''---and this should apply to its usage in propositional logic, as well as in predicate logic---is better abandoned if it is tied to the philosophical position that logic is devoid of content, that it says nothing. As the rest of mathematics, logic should be indifferent towards this position. (A citation covering this matter is in \cite{CN09}, p.\ 73.)

In making this comment, G\" odel might have had in mind the author usually credited for introducing the term ``tautology'' in modern logic; namely, Wittgenstein (in \cite{W21}; Kant, however, already used the term in \cite{K00}, Section I.37, for a particular kind of analytic propositions).
The philosophical opinions of the early Wittgenstein, and perhaps also of the later, as well as those of the logical positivists, could be described by saying that they thought that logic, and mathematics too, are devoid of content.

\vspace{2ex}

\noindent \textbf{\S 14. ``Thinking machines''.} G\" odel speaks of matters of completeness and decidability from \S 11 with a simile (see \cite{CN09}, Extract~2, pp.\ 84-85). Turing machines are not named, but their working is suggested by, in G\" odel's words, ``thinking machines''. One device has a crank, which has to be turned to produce tautologies, where this word covers also the universally true formulae of predicate logic, and there is another device with a typewriter and a bell, which rings if one types in a tautology. The device of the second kind is available for propositional logic and monadic predicate logic, but not for the whole predicate calculus. G\" odel says: ``So here already one can prove that [the] Leibnitzian program of the 'calculemus' cannot be carried
through; i.e., one knows that the human mind will never be able to be replaced by a machine already for this comparatively simple question to decide whether a formula is a tautology or not.''

\vspace{2ex}

\noindent \textbf{\S 15. Existential presuppositions.} G\" odel rejects the unjustified existential presuppositions of Aristotelian syllogistic. This he does because they are either an empirical matter foreign to logic, or they would hamper arguments where they are not made, in which the issue might be exactly whether a predicate applies to anything or not, as it happens in mathematics.

\vspace{2ex}

\noindent \textbf{\S 16. Classes and relations.} G\" odel begins the last part of the course by dealing with matters that lead to the theory of types. He first considers classes, which he does not call \emph{sets}, as extensions of unary predicates, and briefly mentions Russell's ``no class theory'', referred to in \cite{DA16} (Section~4). He defines the usual Boolean operation on classes, pointing out the similarities between these and arithmetical operations. Relations as extensions of predicates of arity greater than one are also introduced, together with some related elementary notions, like symmetry and transitivity, and operations, like inverse and composition, are considered. The notion of function is introduced in this context. (Much of this stuff is crossed out in the notes.) Russell's convention for understanding definite descriptions is presented in a few sentences printed in \cite{CN09} (pp.\ 70-71), and Russell's achievement is succinctly presented as making the meaningfulness of language independent of empirical matters.

In the discussion of Aristotelian syllogistic in G\" odel's notes one finds a formal theory of propositional logic based on a language with the individual variables $\alpha,\beta,\ldots$, which should be interpreted by sets of objects, and two binary relations \textbf{a} and \textbf{i}, such that $\alpha \,{\mathbf a}\,\beta$ should mean that $\alpha$ is a subset of $\beta$, and $\alpha \,{\mathbf i}\,\beta$ should mean that $\alpha$ and $\beta$ have a nonempty intersection. With the axioms stating the reflexivity and transitivity of \textbf{a}, transitivity corresponding to the mood Barbara, and in addition $(\alpha\, {\mathbf i}\,\beta\; .\; \beta\, {\mathbf a}\,\gamma)\supset\gamma\, {\mathbf i}\,\alpha$, which corresponds to the mood Dimatis (or Dimaris), added to the axioms and rules of propositional logic, it is claimed that one obtains a complete system for Aristotelian logic without the restriction to nonempty sets of objects. (The Aristotelian \textbf{e} and \textbf{o} relations are defined in terms of \textbf{i}, \textbf{a} and negation.) It is easy to see that G\" odel's Dimatis axiom, could be replaced by $(\beta\, {\mathbf a}\,\gamma\; .\; \beta\, {\mathbf i}\,\alpha)\supset\alpha\, {\mathbf i}\,\gamma$, which corresponds to the mood Datisi. (One derives from either Dimatis or Datisi and the reflexivity of \textbf{a} that \textbf{i} is symmetric.)

The system with the Datisi axiom and in addition the reflexivity of \textbf{i} may be found in {\L}ukasiewicz's book \cite{L50} (\S 25), whose investigations of these matters date from the summer of 1939, a few months after G\" odel's course in Notre Dame. {\L}ukasiewicz's system axiomatizes Aristotelian logic with the usual restriction to nonempty sets of objects. A proof that {\L}ukasiewicz's system without the reflexivity of \textbf{i} and with the additional axioms $\alpha\, {\mathbf i}\,\beta \supset \alpha\, {\mathbf i}\,\alpha$ and $\alpha\, {\mathbf i}\,\alpha \vee \alpha\, {\mathbf a}\,\beta$ is complete for the interpretation envisaged by G\" odel may be found in Shepherdson's paper \cite{Sh56}.

The attention G\" odel paid to Aristotelian logic understood within predicate logic may also be explained by G\" odel's expecting that one half of his audience would be made of older philosophers of Notre Dame. According to \cite{Daw97} (p.\ 135) and \cite{Daw}, only the other, mathematical, half, of some ten younger graduate students, attended the lectures till the end.

\vspace{2ex}

\noindent \textbf{\S 17. Paradoxes and type theory.} After all that, G\" odel turns to paradoxes, of which he mentions the Burali-Forti paradox of the set of all ordinal numbers and Russell's paradox of the set of all sets that are not members of themselves. The latter is treated in more detail, and there is a brief discussion of the theory of types. G\" odel's remarks on paradoxes are published as Extract~3 in \cite{CN09} (pp.\ 85-89). G\" odel's opinions on self-reference here point to his interest for the intensional logic of concepts, which is a matter we have discussed in \cite{DA16} (Section~5). With that, G\" odel's course ends (if we don't count jottings concerning matters covered before).

\end{document}